\documentclass[12pt]{article}

\usepackage{amsmath,amssymb,amsbsy,amsfonts,amsthm,latexsym,
                    amsopn,amstext,amsxtra,euscript,amscd}

\begin{document}

    \textwidth 5.35 in
   \textheight 8.3 in
\newtheorem{theorem}{Theorem}
\newtheorem{lemma}[theorem]{Lemma}
\newtheorem{claim}[theorem]{Claim}
\newtheorem{cor}[theorem]{Corollary}
\newtheorem{prop}[theorem]{Proposition}
\newtheorem{definition}{Definition}
\newtheorem{question}[theorem]{Question}

\def\cA{{\mathcal A}}
\def\cB{{\mathcal B}}
\def\cC{{\mathcal C}}
\def\cD{{\mathcal D}}
\def\cE{{\mathcal E}}
\def\cF{{\mathcal F}}
\def\cG{{\mathcal G}}
\def\cH{{\mathcal H}}
\def\cI{{\mathcal I}}
\def\cJ{{\mathcal J}}
\def\cK{{\mathcal K}}
\def\cL{{\mathcal L}}
\def\cM{{\mathcal M}}
\def\cN{{\mathcal N}}
\def\cO{{\mathcal O}}
\def\cP{{\mathcal P}}
\def\cQ{{\mathcal Q}}
\def\cR{{\mathcal R}}
\def\cS{{\mathcal S}}
\def\cT{{\mathcal T}}
\def\cU{{\mathcal U}}
\def\cV{{\mathcal V}}
\def\cW{{\mathcal W}}
\def\cX{{\mathcal X}}
\def\cY{{\mathcal Y}}
\def\cZ{{\mathcal Z}}

\def\A{{\mathbb A}}
\def\B{{\mathbb B}}
\def\C{{\mathbb C}}
\def\D{{\mathbb D}}
\def\E{{\mathbb E}}
\def\F{{\mathbb F}}
\def\G{{\mathbb G}}
\def\I{{\mathbb I}}
\def\J{{\mathbb J}}
\def\K{{\mathbb K}}
\def\L{{\mathbb L}}
\def\M{{\mathbb M}}
\def\N{{\mathbb N}}
\def\O{{\mathbb O}}
\def\P{{\mathbb P}}
\def\Q{{\mathbb Q}}
\def\R{{\mathbb R}}
\def\S{{\mathbb S}}
\def\T{{\mathbb T}}
\def\U{{\mathbb U}}
\def\V{{\mathbb V}}
\def\W{{\mathbb W}}
\def\X{{\mathbb X}}
\def\Y{{\mathbb Y}}
\def\Z{{\mathbb Z}}

\def\ep{{\mathbf{e}}_p}
\def\em{{\mathbf{e}}_m}

\def\scr{\scriptstyle}
\def\\{\cr}
\def\({\left(}
\def\){\right)}
\def\[{\left[}
\def\]{\right]}
\def\<{\langle}
\def\>{\rangle}
\def\fl#1{\left\lfloor#1\right\rfloor}
\def\rf#1{\left\lceil#1\right\rceil}
\def\le{\leqslant}
\def\ge{\geqslant}
\def\eps{\varepsilon}
\def\mand{\qquad\mbox{and}\qquad}

\def\vec#1{\mathbf{#1}}
\def\inv#1{\overline{#1}}
\def\vol#1{\mathrm{vol}\,{#1}}
\def\dist{\mathrm{dist}}

\def\SL{\mathrm{SL}}

\def\Hba{\overline{\cH}_{a,m}}
\def\Hta{\widetilde{\cH}_{a,m}}
\def\Hb1{\overline{\cH}_{m}}
\def\Ht1{\widetilde{\cH}_{m}}

\def\Zm{\Z/m\Z}

\def \tX{\widetilde{X}}
\def\tY{\widetilde{Y}}
\def\tZ{\widetilde{Z}}

\def\Err{{\mathbf{E}}}

\newcommand{\comm}[1]{\marginpar{%
\vskip-\baselineskip 
\raggedright\footnotesize
\itshape\hrule\smallskip#1\par\smallskip\hrule}}

\def\xxx{\vskip5pt\hrule\vskip5pt}


\title{\bf On the Restricted Divisor Function in 
Arithmetic Progressions}

\author{
{\sc Igor E. Shparlinski} \\
{Department of Computing, Macquarie University} \\
{Sydney, NSW 2109, Australia} \\
{igor.shparlinski@mq.edu.au}}

\date{\today}
\pagenumbering{arabic}

\maketitle

\begin{abstract}
We obtain several asymptotic estimates for the sums 
of the restricted divisor function
$$
\tau_{M,N}(k) = \# \{1 \le m \le M, \ 1\le n \le N: mn = k\}
$$ 
over short arithmetic progressions,  which improve some
results of J.~Truelsen. Such estimates are motivated by 
the links with the pair correlation problem for fractional 
parts of the quadratic function $\alpha k^2$, $k=1,2,\ldots$ 
with a real $\alpha$. 
\end{abstract}



\section{Introduction}
\label{sec:intro}

There is a long history of studying the distribution of the
divisor function over short arithmetic progressions,
see~\cite{BHS,Blom,Fouv,FrIw1,FrIw2} and references therein.

Recently, Truelsen~\cite{True} has introduced the restricted 
divisor function 
$$
\tau_{M,N}(k) = \# \{1 \le m \le M, \ 1\le n \le N: mn = k\}
$$
and shown its relevance to the {\it pair correlation problem\/} for 
fractional parts of  the quadratic
function $\alpha k^2$, $k=1,2,\ldots$, with a real $\alpha$,
see also~\cite{HB,RudSar} for various results  and conjectures 
concerning this problem. In particular, it is conjectured
in~\cite[Conjecture~1.2]{True} that for any fixed  
$\varepsilon, \delta, c_1, c_2>0$, 
if positive integers $N,M,R$ and $q$ satisfy
$$
 N \ge q^{1/2+\varepsilon}
 \qquad c_1N \le M \le c_2 N  \qquad R \ge N^\delta,
$$
then, uniformly over all integers $a$ with $\gcd(a,q)=1$, we have
\begin{equation}
\label{eq:Conj}
\sum_{r=1}^R \sum_{k \equiv ar \pmod q} \tau_{M,N}(k) 
\sim \frac{MNR}{q}.
\end{equation}
It is also shown in~\cite{True} that the asymptotic 
formula~\eqref{eq:Conj} yields explicit examples 
of real $\alpha$ for which distribution of spacings between 
the 
fractional parts of  $\alpha k^2$ is {\it Poissonian\/}.  

Towards the conjecture~\eqref{eq:Conj}, several asymptotic
formulas and estimates are derived in~\cite{True}.

In particular, as in~\cite{True}, for positive integer  $q$, $M$, $N$  and
a divisor $d\mid q$,    we consider the sums
\begin{equation}
\label{eq:Delta}
\Delta_q(d;M,N) = \sum_{\substack{a=1\\\gcd(a,q)=d}}^q
\left| \sum_{k \equiv a\pmod q} \tau_{M,N}(k)- \frac{MN}{q^2} \Phi(q,d)\right|^2,
\end{equation}
where 
\begin{equation}
\label{eq:Phi}
\Phi(q,d) = \sum_{e\mid d} e
 \sum_{f|q/e} f \mu\(\frac{q}{ef}\) =
 \sum_{e\mid d} e \varphi(q/e),
\end{equation}
see~\cite[Equation~(1.5)]{True} and  $\mu(k)$ is the M{\"o}bius function.
Also as in~\cite{True}, for positive integer  $q$, $M$, $N$ and $R$, we  consider the sums 
\begin{equation}
\label{eq:Gamma}
\Gamma_q(M,N,R) = 
\sum_{a=1}^{q-1}
\left|
\sum_{r\le R}\sum_{k \equiv ar\pmod q} \tau_{M,N}(k)- \frac{MNR}{q}\right|^2.
\end{equation}

Here, in  Section~\ref{sec:Aver1}, we show that a result of~\cite{Shp} almost instantly 
implies the estimate of~\cite[Theorem~1.8]{True} on 
the sums $\Delta_q(d;M,N)$, and in fact, 
in a slightly stronger form. Furthermore, using a
different technique of multiplicative character sums, 
in  Section~\ref{sec:Aver2}
we obtain a new estimates on the sums $\Gamma_q(M,N,R)$, which for some parameter ranges
improves that of~\cite[Theorem~1.9]{True}. We present our argument 
only in the case of prime $q$ but combining it with elementary
(but somewhat cluttered) sieving it can also be used for arbitrary 
$q$.

\section{Preliminaries}

\subsection{General notation and facts} 

Throughout the paper, any implied constants in symbols $O$, $\ll$
and $\gg$ may occasionally depend on the positive parameters $\varepsilon$ 
and $\delta$ and  are absolute otherwise. We recall
that the notations $U = O(V)$,  $U \ll V$ and  $V \gg U$  are
all equivalent to the statement that $|U| \le c V$ holds 
with some constant $c> 0$.

We always assume that the variables which appear in 
congruences and as arguments of standard arithmetic 
functions are  integers.

We recall that for 
$$
\varphi(s,K)  = \sum_{\substack{1 \le k \le K\\ \gcd(k,s)=1}} 1
$$
we have the asymptotic formula
\begin{equation}
\label{eq:phiqM}
\varphi(s,K) 
= \frac{\varphi (s) }{s}K + O(s^{o(1)}), 
\end{equation}
see~\cite[Equation~(3.1)]{True}, that  follows  from the 
inclusion-exclusion principle and the 
well-known bound on the divisor and Euler functions
$$
\tau(s) = s^{o(1)} \mand \varphi(s) = s^{1+o(1)},
$$
see~\cite[Theorem~317]{HardyWright} 
and~\cite[Theorem~328]{HardyWright},
respectively. 

\subsection{Character sums}
\label{sec:Char}

Let $\varPhi_s$ be the set of all $\varphi(s)$
multiplicative characters modulo  $s$. 
We also use $\chi_0$ to denote the principal character
and  
$$
\varPhi_s^* = \varPhi_s \setminus \{\chi_0\}
$$
to denote the set of nonprincipal multiplicative 
characters modulo $s$.

For an integer $Z$ and $\chi \in \varPhi_s$ we define the sums
\begin{equation}
\label{eq:Sum Sst}
S_s(Z;\chi) = \sum_{z = 1}^{Z} \chi(z).
\end{equation}

The following result  is a combination of  the
P{\'o}lya-Vinogradov  (for $\nu =1$) and Burgess
(for $\nu\ge2$)  bounds,
see~\cite[Theorems~12.5 and 12.6]{IwKow}.

\begin{lemma}
   \label{lem:PVB} For a prime $s$ and  positive integers $Z \le s$, 
   the bound
$$
\max_{\chi \in \varPhi_s^*}
\left| S_s(Z;\chi)\right|  \le Z^{1 -1/\nu} s^{(\nu+1)/4\nu^2 + o(1)}
$$
holds with an arbitrary fixed 
integer $\nu\ge 1$.
\end{lemma}

We combine Lemma~\ref{lem:PVB} with a bound on 
the fourth moment of the sums $S_s(Z,t;\chi)$.
First we recall the
following   estimate from~\cite{ACZ}
(for prime $s$) and~\cite{FrIw2} 
(for arbitrary $s$), see also~\cite{CochSih,GarGar}, 
which we present in the 
following slightly  relaxed form.

\begin{lemma}
   \label{lem:4th Moment} 
   For positive integers $Z \le s$, the bound
$$
\sum_{\chi \in \varPhi_s^*}
\left| S_s(Z;\chi)\right|^4
\le   s^{1 + o(1)} Z^2
$$
holds.
\end{lemma}

\subsection{Sums with $\tau_{M,N}(k)$ and congruences}
\label{sec:prod res}

We note that sums of the restricted divisor function 
over an arithmetic progression can be expressed via the number 
of solutions to a certain congruence. For example,
\begin{equation}
\label{eq:T tau}
\sum_{k \equiv a \pmod q} \tau_{M,N}(k) 
= T_q(M,N;a), 
\end{equation}
where $T_q(M,N;a)$ is number of solutions to the congruence
\begin{equation}
\label{eq:Cong mn}
mn\equiv a \pmod q, \qquad 1 \le m \le M, \ 1 \le n \le N.
\end{equation}
This interpretation underlines our approach.

To estimate the function $T_s(M,N;a)$  
it is more convenient to work with the  quantity 
$T_s^*(X,Y;a)$ which is defined as number of solutions to the congruence
$$
xy\equiv a \pmod s, \qquad 1 \le x \le X,\ \gcd(x,s)=1, \ 1 \le y \le Y.
$$

One of our main tool is the following
special case of~\cite[Theorem~1]{Shp}, 
combined with~\eqref{eq:phiqM}.

\begin{lemma}
\label{lem:T Aver}  For  positive integers $s$, $X\le Y$, 
we have 
$$
\sum_{a=1}^s \left| T_s^*(X,Y;a) -
 \frac{\varphi(s)}{s^2}  XY \right|^2
\le  XY  s^{o(1)}.
$$
\end{lemma}

We also define $R_s(X,Y,Z;a)$ 
as number of solutions to the congruence
$$
xy\equiv az \pmod{s}, 
$$
with 
$$ 1 \le x \le X, \qquad 1 \le y \le Y, \qquad 
\ 1 \le z \le Z.
$$ 

\begin{lemma}
\label{lem:R Aver}  For a prime $s$ and positive integers 
 $X,Y,Z < s$ 
we have 
$$
\sum_{a=1}^{s-1} \left|R_s(X,Y,Z;a) -
\frac{XYZ}{s-1}\right|^2 
\le XYZU^{1 -2/\nu} s^{(\nu+1)/2\nu^2+o(1)}   
$$
where  $U = \min\{X,Y,Z\}$  and  $\nu\ge 1$ is arbitrary  fixed 
positive integer.
\end{lemma}

\begin{proof}
We note that for every $a$ with $\gcd(a,s)=1$, we obtain
$$
R_s(X,Y,Z;a)=
\frac{1}{s-1} 
\sum_{x=1}^X
\sum_{y=1}^Y \sum_{z=1}^Z
   \sum_{\chi \in \varPhi_s}\chi\(a^{-1}xyz^{-1}\).
$$

Recalling the definition~\eqref{eq:Sum Sst}, 
changing the order of summation, using that 
$$
\chi\(z^{-1}\) = \overline{\chi}(z), 
$$
if $\gcd(z,s)=1$ 
where $\overline{\chi}$ is the complex conjugated character,
we derive
$$
R_s(X,Y,Z;a) =
\frac{1}{s-1} 
\sum_{\chi \in \varPhi_s}\overline{\chi}\(a\) S_s(X;\chi) 
S_s(Y;\chi) S_s(Z;\overline{\chi}) .
$$
We now  separate the contribution from the 
principal character $\chi =\chi_0$, getting
\begin{equation*}
\begin{split}
R_s(X,Y,Z;a)-  
\frac{XYZ}{s-1}
= \frac{1}{s-1} 
\sum_{\chi \in \varPhi_s*}\overline{\chi}\(a\)  S_s(X;\chi) 
S_s(Y;\chi) S_s(Z;\overline{\chi}) .
\end{split}
\end{equation*}
 
Using the orthogonality of characters, we easily derive
\begin{equation*}
\begin{split}
\sum_{a=1}^s \left| R_s(X,Y,Z;a) -
\frac{XYZ}{s-1}\right|^2&\\
 = \frac{1}{\varphi(s)} 
\sum_{\chi  \in \varPhi_s*} &| S_s(X,u;\chi)|^2 
|S_s(Y,v;\chi)|^2 |S_s(Z;w,\chi)|^2\\
 = \frac{1}{\varphi(s)} 
\sum_{\chi  \in \varPhi_s*} &| S_s(\tX,u;\chi)|^2 
|S_s(\tY,v;\chi)|^2 |S_s(\tZ;w,\chi)|^2, 
\end{split}
\end{equation*}
for any  permutation $(\tX,\tY,\tZ)$ is any of $(X,Y,Z)$.
We now apply Lemma~\ref{lem:PVB} to the last sum and then
use  the Cauchy inequality, arriving to
\begin{equation*}
\begin{split}
\sum_{a=1}^s  \left| R_s(X,Y,Z;a) -
\frac{XYZ}{s-1}\right|^2&\\
 \le  \frac{\tZ^{2 -2/\nu} s^{(\nu+1)/2\nu^2+o(1)}}{s-1} &
\sqrt{\sum_{\chi  \in \varPhi_s*}|S_s(\tX;\chi)|^4}
\sqrt{\sum_{\chi  \in \varPhi_s*}|S_s(\tY;\chi)|^4}.
\end{split}
\end{equation*}
We now choose a permutation $(\tX,\tY,\tZ)$ with  
$\tZ =  U = \min\{X,Y,Z\}$. 
Using Lemma~\ref{lem:4th Moment}, we obtain the desired result.
\end{proof}

\section{Average Values $\tau_{M,N}(k)$ over Some Families of Progressions}
\label{sec:Aver}

\subsection{One parameter family of progressions}
\label{sec:Aver1}

Here we estimate the sums $\Delta_q(d;M,N)$ given by~\eqref{eq:Delta}. 
and  show how Lemma~\ref{lem:T Aver} implies
a stronger and more general form of 
the estimate~\cite[Theorem~1.8]{True} which asserts 
that if $M \ll N \ll M$ then 
\begin{equation}
\label{eq:True bound}
\Delta_q(d;M,N) 
\le \frac{1}{q} N^{\max\{7/2, 4-\delta\} + o(1)},
\end{equation}
uniformly over $q \le N^{2-\delta}$ and $d\mid q$

\begin{theorem}
\label{thm:Delta} For arbitrary positive integers $q$, $M$ and
$N$ and a divisor $d\mid q$ we have 
$$
\Delta_q(d;M,N) 
\le MN q^{o(1)}
$$
\end{theorem}

\begin{proof} 
Without loss of generality we can assume that 
$M \ge N$.

For each  divisor $e\mid d$, we collect together
the solutions to~\eqref{eq:Cong mn} with $\gcd(m,q)=e$, 
getting 
$$
T_q(M,N;a) = \sum_{e\mid d} T_{q/e}^*(\fl{M/e},N;a/e) .
$$
where $T_s^*(X,Y;a)$ is defined in Section~\ref{sec:prod res}.

Recalling~\eqref{eq:T tau}  and~\eqref{eq:Phi}, we obtain
\begin{eqnarray*}
\Delta_q(d;M,N) & = &
\sum_{\substack{a=1\\\gcd(a,q)=d}}^q
\left| T_q(M,N;a) - \frac{MN}{q^2} \Phi(q,d)\right|^2\\
& = &
\sum_{\substack{a=1\\\gcd(a,q)=d}}^q
\(
\sum_{e\mid d} \left|T_{q/e}^*(\fl{M/e},N;a/e) -  \frac{MNe}{q^2}  \varphi(q/e)\right|\)^2.
\end{eqnarray*}

Thus, using the Cauchy inequality, we obtain
\begin{equation}
\label{eq:Prelim}
\begin{split}
\Delta_q(d;M,N&)
 \\
\le
&q^{o(1)} \sum_{e\mid d}
\sum_{\substack{a=1\\\gcd(a,q)=d}}^q
\left| T_{q/e}^*(\fl{M/e},N;a/e) 
-  \frac{MNe}{q^2}  \varphi(q/e)\right|^2.
\end{split}
\end{equation}

We now note that
\begin{eqnarray*}
\lefteqn{\frac{MNe}{q^2}\varphi(q/e)  =\frac{(M/e)N}{(q/e)^2} \varphi(q/e)} \\
 & & \qquad =
  \frac{\fl{M/e}N}{(q/e)^2} \varphi(q/e) + O(Ne/q)
  = \frac{\fl{M/e}N}{(q/e)^2} \varphi(q/e) + O(Nd/q).
\end{eqnarray*} 
We now see from~\eqref{eq:Prelim} that
\begin{equation*}
\begin{split}
\Delta_q(d;M,&N)
 \\
\le
&q^{o(1)} \sum_{e\mid d}
\sum_{\substack{a=1\\\gcd(a,q)=d}}^q
\left|T_{q/e}^*(\fl{M/e},N;a/e) 
- \frac{\fl{M/e}N}{(q/e)^2} \varphi(q/e)\right|^2 \\
& \qquad\qquad\qquad\qquad\qquad\qquad\qquad\qquad + N^2 d^2 q^{-2+o(1)} \sum_{e\mid d}
\sum_{\substack{a=1\\\gcd(a,q)=d}}^q 1\\
\le
&q^{o(1)} \sum_{e\mid d}
\sum_{\substack{a=1\\\gcd(a,q)=d}}^q 
\left| T_{q/e}^*(\fl{M/e},N;a/e) 
- \frac{\fl{M/e}N}{(q/e)^2} \varphi(q/e)\right|^2 \\
& \qquad\qquad\qquad\qquad\qquad\qquad\qquad\qquad\qquad\qquad
 + N^2d q^{-1+o(1)}.
 \end{split}
\end{equation*}
Writing $a=ce$, we derive
\begin{equation*}
\begin{split}
\Delta_q(d;M,N)&\\
\le &q^{o(1)} \sum_{e\mid d}
\sum_{c=1}^{q/e} 
\left|T_{q/e}^*(\fl{M/e},N;c) 
- \frac{\fl{M/e}N}{(q/e)^2} \varphi(q/e)\right|^2 \\
& \qquad\qquad\qquad\qquad\qquad\qquad\qquad\qquad \qquad\qquad+ N^2d q^{-1+o(1)}  .
\end{split}
\end{equation*}
Now  recalling Lemma~\ref{lem:T Aver}, 
we obtain
$$
\Delta_q(d;M,N) \le MNq^{o(1)} + N^2d q^{-1+o(1)} \le (M+N)Nq^{o(1)} .
$$
Since $M \ge N$, this concludes the proof. 
\end{proof}

Note  that  the bound of Theorem~\ref{thm:Delta}
is more general than~\eqref{eq:True bound} as it
works for $M$ and $N$ of essentially different sizes and
does not need the restriction. In particular,
if $M \ll N \ll M$, 
then this  bound takes form $N^{2} q^{o(1)}$, which 
improves~\eqref{eq:True bound} for $\delta > 1/2$, that 
is, for $N \ge q^{2/3+\varepsilon}$ for any fixed
$\varepsilon > 0$.

\subsection{Two parameter family of progressions}
\label{sec:Aver2}

Note that in~\cite{True} the bound~\eqref{eq:True bound} 
has been used to prove several other results. 
Theorem~\ref{thm:Delta} can be used to get corresponding
generalisations and improvements of these bounds. For example,
bounds of $\Delta_q(d;M,N)$ are used 
in~\cite[Theorem~1.9]{True} to derive the estimate
on the sums  $\Gamma_q(M,N,R)$ given by~\eqref{eq:Gamma}. In 
particular, by~\cite[Theorem~1.9]{True} we have
\begin{equation}
\label{eq:True-Gamma}
\Gamma_q(M,N,R)  \le  N^4 R^2
\(R^{-2} + N^{\max\{-1/2,-\delta\}}\) q^{-1+o(1)}, 
\end{equation}
provided $M \ll N \ll M$, $R \le q \le N^{2-\delta}$
(note that the condition of~\cite[Theorem~1.9]{True} 
that $R\ge N^\eta$ for some positive $\eta>0$
does not seem to be needed for the bound, but the bound is 
nontrivial only if it is satisfied). 
The estimate~\eqref{eq:True-Gamma} shows
that the conjectured asymptotic formula~\eqref{eq:Conj} holds 
on average under appropriate averaging conditions, 
see~\cite[Corollary~1.10]{True}. 

As in the case of $\Delta_q(d;M,N)$, using Theorem~\ref{thm:Delta} 
one now obtains a similar 
generalisation and improvement for $\Gamma_q(M,N,R)$. 
One can probably use similar arguments to  
sharpen~\cite[Theorem~4.5]{True} as well.

Furthermore, we now present a different approach,
based on Lemma~\ref{lem:R Aver}, which allows us to
obtain estimates on $\Gamma_q(M,N,R)$ that  are sometimes
stronger that those of~\cite[Theorem~1.9]{True}  or 
following from Theorem~\ref{thm:Delta}. We demonstrate this
approach only in the case of prime modulus $q$. In the general 
case, one can use it as well, but it involves rather 
cluttered expressions arising from the inclusion-exclusion
principle.

\begin{theorem}
\label{thm:Gamma} 
For a prime $q$ and positive integers 
$M,N,R < q$, the bound 
$$
\Gamma_q(M,N,R) \le MNR L^{1 -2/\nu}q^{(\nu+1)/2\nu^2+o(1)}  
$$
holds, where  $L = \min\{M,N,R\}$  and  $\nu\ge 1$ is arbitrary  fixed 
positive integer.

\end{theorem}

\begin{proof} As in the proof of Theorem~\ref{thm:Delta} we see
that 
$$
\Gamma_q(M,N,R) =
\sum_{a=1}^{q-1}
\left| R_q(M,N,R;a) - \frac{MNR}{q} \right|^2, 
$$
and using Lemma~\ref{lem:R Aver},  we conclude the proof. 
\end{proof}

For example, if $q$ is prime then for $M,N= q^{2/3 + o(1)}$ 
and $R =q^{1/2 + o(1)}$, 
applying Theorem~\ref{thm:Gamma} with $\nu = 2$ we obtain 
$$
\Gamma_q(M,N,R) \le  q^{53/24+o(1)} 
$$
while~\eqref{eq:True-Gamma} gives only 
$$
\Gamma_q(M,N,R) \le  q^{7/3+o(1)} 
$$
for the above choice of parameters. One can certainly easily produce
many other examples of the parameters $(M,N,R)$ for which 
 Theorem~\ref{thm:Gamma}  is stronger than~\eqref{eq:True-Gamma}.
 
As we have said the argument used in the proof of Lemma~\ref{lem:R Aver}
an thus of  Theorem~\ref{thm:Gamma} can also be applied in the case 
of composite $q$. However we recall that the Burgess bound for character sums modulo a composite $q$ has some 
limitations on the possible choices of $\nu$, see~\cite[Theorem~12.6]{IwKow} for details.

\section*{Acknowledgement}

The author is grateful to Jimi Truelsen and to the referee for 
very useful comments.  During the preparation of this paper, the author
was supported in part by ARC grant DP1092835.

\end{document}